# TOROIDAL COMPACTIFICATIONS AND STACKY COHOMOLOGY OF MUMFORD–TATE DOMAINS

MOHAMMAD REZA RAHMATI

Abstract. Mumford–Tate domains parametrize polarized Hodge structures with fixed Mumford–Tate group and play a central role in the geometry of period maps. Their degenerations are governed by nilpotent orbits and limiting mixed Hodge structures, whose asymptotics are encoded in the logarithmic compactifications of Kato–Usui. In this paper we construct and study the *log–toric Hodge stack*

$$\mathcal{D}^{\log}_{M,\Sigma} := [D_{M,\Sigma}/\Gamma],$$

obtained from a Mumford–Tate domain $D_M$ and a fan $\Sigma$ of nilpotent cones by forming the quotient of the Kato–Usui partial compactification $D_{M,\Sigma}$ by a neat arithmetic group $\Gamma \subset M(\mathbb{Q})$. We show that $\mathcal{D}^{\log}_{M,\Sigma}$ is a global quotient Deligne–Mumford stack, that it admits a natural logarithmic structure extending the period domain, and that near every boundary stratum associated to a cone $\sigma \in \Sigma$ it admits a canonical analytic log–étale chart of the form

$$\Big([F_\sigma/G_\sigma] \times \mathcal{T}_\sigma\Big)^\circ,$$

where $F_\sigma$ is the space of nilpotent orbits modulo unipotent actions, $G_\sigma$ is a finite symmetry group of the associated limiting mixed Hodge structures, and $\mathcal{T}_\sigma$ is a toric Deligne–Mumford stack refining the toric variety $D_\sigma$ attached to $\sigma$. This decomposition cleanly separates Hodge-theoretic information from the combinatorial and stacky boundary data.

The log–toric structure yields a natural home for characteristic cohomology. Using Behrend's comparison theorem for differential forms on quotient stacks, we identify the de Rham complex of $\mathcal{D}^{\log}_{M,\Sigma}$ with the $\Gamma$–equivariant de Rham complex of $D_{M,\Sigma}$. We introduce logarithmic Hodge bundles over $\mathcal{D}^{\log}_{M,\Sigma}$ and formulate a stacky form of the Green–Griffiths–Kerr characteristic cohomology, showing that its classes arise from Chern classes of these bundles together with cohomology classes of toric boundary strata. This leads to a natural *Hodge–generation conjecture* predicting that $H^\bullet(\mathcal{D}^{\log}_{M,\Sigma}, \mathbb{Q})$ is generated by these classes with all relations coming from universal Chern identities and toric stack relations.

Finally, we identify the stack-theoretic stabilizers of boundary points with finite subgroups of the Mumford–Tate groups of the corresponding limiting mixed Hodge structures, commuting with monodromy. This description recovers classical boundary symmetries in rank one (e.g. roots of unity acting on the Tate curve) and provides a conceptual framework for the finite automorphism groups that arise in higher-rank degenerations. Together, these results exhibit $\mathcal{D}^{\log}_{M,\Sigma}$ as a natural global object encoding both the toroidal geometry and the Hodge-theoretic structure of Mumford–Tate domains, providing a unified approach to their characteristic cohomology, degenerations, and boundary symmetries.

---

## 1. Introduction

The global geometry of period maps and their degenerations lies at the heart of modern Hodge theory. On the one hand, Griffiths' period domains [13] provide a homogeneous complex-analytic parameter space for polarized Hodge structures of fixed type. On the other hand, degenerations of such structures—and the limiting mixed Hodge structures (LMHS) that control their asymptotics—are governed by nilpotent orbits and logarithmic geometry in the sense of Kato–Usui [16, 15]. The present paper proposes a framework that combines these two aspects into a single geometric object: a *log–toric Hodge stack*

$$\mathcal{D}^{\log}_{M,\Sigma} := [D_{M,\Sigma}/\Gamma]$$

associated to a Mumford–Tate domain and a Kato–Usui fan of nilpotent cones. We show that this stack is locally (in the analytic log–étale topology) a product of a period–type stack and a toric Deligne–Mumford stack, and we use this structure to organize cohomology classes and automorphism groups arising from degenerations.

Let $(H_{\mathbb{Z}}, Q)$ be a polarized integral Hodge structure of weight $w$ and fixed Hodge numbers, and let

$$G := \mathrm{Aut}(H_{\mathbb{Q}}, Q)$$

be the $\mathbb{Q}$–algebraic group of $Q$–preserving automorphisms. Griffiths' period domain $D$ is a homogeneous complex manifold

$$D \;\cong\; G(\mathbb{R})^+/H,$$

parametrizing Hodge filtrations $F^\bullet \subset H_{\mathbb{C}}$ satisfying the Hodge–Riemann bilinear relations [13]. The Hodge tensors of a reference Hodge structure cut out a smaller algebraic subgroup

$$M \subset G,$$

the *Mumford–Tate group*. The corresponding *Mumford–Tate domain*

$$D_M \subset D$$

is again homogeneous,

$$D_M \;\cong\; M(\mathbb{R})^+/H_M,$$

and serves as the natural parameter space for Hodge structures with prescribed Mumford–Tate group in the sense of Green–Griffiths–Kerr [11, 12] and Han–Robles [14].

Fix a neat arithmetic subgroup $\Gamma \subset M(\mathbb{Q})$, so that the quotient

$$S^\circ := \Gamma\backslash D_M$$

is a smooth complex-analytic space. A key theme in recent work is that the geometry and arithmetic of $S^\circ$ become much more tractable when one passes to the level of stacks and when one incorporates a canonical compactification: stack-theoretic techniques clarify cohomology and intersection theory [6, 5, 2, 1], while partial compactifications in the sense of Kato–Usui provide a natural home for LMHS and boundary behavior [16, 15, 17].

The Kato–Usui theory attaches to $(D_M, \Gamma)$ and a fan $\Sigma$ of nilpotent cones in $\mathrm{Lie}(M)_{\mathbb{Q}}$ a logarithmic analytic space

$$D_{M,\Sigma},$$

containing $D_M$ as a dense open subset and equipped, for each cone $\sigma \in \Sigma$, with a log–smooth map

$$\pi_\sigma : U_\sigma \longrightarrow D_\sigma$$

to an affine toric variety $D_\sigma$ built from the monodromy monoid $\Gamma(\sigma) \subset \Gamma$ [16]. Nilpotent orbits and LMHS with monodromy in $\sigma$ are parametrized over $D_\sigma$, and the various $U_\sigma$ glue to a partial compactification $D_{M,\Sigma}$ whose boundary reflects the combinatorics of the fan $\Sigma$.

From the perspective of stacks, it is natural to pass from the spaces $D_M$ and $D_{M,\Sigma}$ to the quotient stacks

$$\mathcal{D}_M := [D_M/\Gamma], \qquad \mathcal{D}^{\log}_{M,\Sigma} := [D_{M,\Sigma}/\Gamma].$$

These are analytic Deligne–Mumford stacks with coarse moduli spaces $S^\circ$ and $S_\Sigma := \Gamma\backslash D_{M,\Sigma}$, respectively. Behrend's theory identifies their de Rham cohomology with equivariant cohomology [6, 5], providing a natural home for characteristic classes of Hodge bundles and for Green–Griffiths–Kerr universal characteristic cohomology [11]. At the same time, the boundary pieces $D_\sigma$ are toric, and their natural stacky refinements are toric Deligne–Mumford stacks in the sense of Geraschenko–Satriano and Fantechi–Mann–Nironi [9, 10, 8], whose cohomology rings and intersection theory are explicitly computable.

The main structural observation of this paper is that these ingredients fit together very cleanly: locally along the boundary, the log–toric Hodge stack $\mathcal{D}^{\log}_{M,\Sigma}$ is a *product* of a period–type stack and a toric stack. More precisely, for each cone $\sigma \in \Sigma$ there is a manifold $F_\sigma$ parametrizing nilpotent orbits modulo unipotent action and a toric Deligne–Mumford stack $\mathcal{T}_\sigma$ refining $D_\sigma$, such that

$$\mathcal{D}^{\log}_{M,\Sigma} \quad \text{is locally isomorphic to} \quad \big([F_\sigma/G_\sigma] \times \mathcal{T}_\sigma\big)^\circ$$

in the analytic log–étale topology, for some finite group $G_\sigma$. All of the stackiness along the boundary is encoded by the toric stacks $\mathcal{T}_\sigma$, while the period–type factor $[F_\sigma/G_\sigma]$ carries the LMHS and Mumford–Tate data.

The first main goal of the paper (accomplished in Section 2) is to lay out this structure systematically. After recalling the Mumford–Tate set-up and the Kato–Usui construction, we define the Hodge stack $\mathcal{D}_M$ and its logarithmic extension $\mathcal{D}^{\log}_{M,\Sigma}$, prove that they are global quotient Deligne–Mumford stacks, and then refine the Kato–Usui local model $U_\sigma \cong F_\sigma \times D_\sigma$ by replacing $D_\sigma$ with a toric Deligne–Mumford stack $\mathcal{T}_\sigma$. The outcome is a precise local log–toric chart

$$\mathcal{U}_\sigma \;\cong\; \big([F_\sigma/G_\sigma] \times \mathcal{T}_\sigma\big)^\circ$$

around each boundary stratum. We illustrate the construction in the classical case of modular curves, where the boundary stack $\mathcal{T}_\sigma$ is a root stack along the cusp, and briefly indicate how higher-rank examples (Siegel modular varieties, orthogonal Shimura varieties) fit into this framework.

The second main goal (Section 3) is to exploit this local product structure to organize two kinds of data that appear ubiquitously in Hodge theory:

- the algebra of cohomology classes generated by Chern classes of Hodge bundles and by toric boundary classes;
- the finite automorphism groups of degenerations (LMHS), viewed as stabilizers of boundary points in the Hodge stack.

On the cohomological side, we introduce logarithmic Hodge bundles $\mathcal{E}^{p,\log}$ on $\mathcal{D}^{\log}_{M,\Sigma}$ extending the usual Hodge bundles on $\mathcal{D}_M$ [16, 15, 17], and we formulate a *stacky Hodge–generation conjecture* stating that the rational cohomology ring $H^\bullet(\mathcal{D}^{\log}_{M,\Sigma}, \mathbb{Q})$ is generated by their Chern classes together with classes of the toric boundary strata, with all relations coming from universal Chern identities and the known Stanley–Reisner and stacky relations in the cohomology of the toric boundary stacks $\mathcal{T}_\sigma$ [9, 8]. Using Behrend's comparison theorem [6] we also give a stacky interpretation of the universal characteristic cohomology of Green–Griffiths–Kerr [11], realized as a natural subcomplex in the de Rham complex of $\mathcal{D}^{\log}_{M,\Sigma}$.

On the group-theoretic side, we identify for each boundary point $x \in \mathcal{D}^{\log}_{M,\Sigma}$ the stacky stabilizer $G_x$ with a finite subgroup of the Mumford–Tate group of the associated LMHS, commuting with the monodromy logarithm $N$:

$$G_x \ \subset \ M_{\mathrm{LMHS}}(\mathbb{Q}) \cap Z(N).$$

In concrete terms, this means that the finite automorphism groups of degenerations are recorded functorially by the stack structure on the boundary. In the classical modular case these are roots of unity acting on the Tate curve; in higher rank they reflect parabolic ranks and torsion in semi-abelian parts of the limit. This perspective complements recent work of Bakker–Brunebarbe–Tsimerman and Bakker–Brunebarbe–Klingler–Tsimerman on definability and GAGA for period maps [3, 4] and is compatible with derived and logarithmic refinements of period maps [7, 1].

Taken together, these results suggest that the log–toric Hodge stack $\mathcal{D}^{\log}_{M,\Sigma}$ provides a natural global home for "Hodge-theoretic symmetries at infinity". Its local decomposition into period and toric factors allows one to separate Hodge–theoretic information from combinatorial toric data, to compute cohomology via Mayer–Vietoris with toric boundary pieces, and to track automorphism groups of degenerations in a stack-theoretic way. In future work we expect this framework to interact fruitfully with the o-minimal and derived approaches to period maps, with the theory of tame and logarithmic stacks [2, 1], and with arithmetic applications involving special subvarieties and unlikely intersections on Mumford–Tate quotients.

## 2. Log–Toric Hodge Stacks: Local and Global Structure

In this section we make precise the idea that a Kato–Usui partial compactification of a Mumford–Tate domain naturally carries the structure of a *log–toric Hodge stack*. Globally, this is a quotient stack

$$\mathcal{D}^{\log}_{M,\Sigma} := [D_{M,\Sigma}/\Gamma]$$

obtained from the Kato–Usui logarithmic space $D_{M,\Sigma}$ and an arithmetic subgroup $\Gamma \subset M(\mathbb{Q})$; locally along the boundary, $\mathcal{D}^{\log}_{M,\Sigma}$ is a product of a period–type stack and a toric Deligne–Mumford stack. We also present a basic example (modular curves) to illustrate the construction.

Throughout we work in the analytic category, following [13, 16, 15, 11, 12]. All stacks are analytic Deligne–Mumford stacks.

### 2.1. Period domains, Mumford–Tate domains, and arithmetic quotients.

Let $(H_\mathbb{Z}, Q)$ be a polarized integral Hodge structure of weight $w$ and fixed Hodge numbers. Denote by $H_\mathbb{Q} := H_\mathbb{Z} \otimes_\mathbb{Z} \mathbb{Q}$ and $H_\mathbb{C} := H_\mathbb{Z} \otimes_\mathbb{Z} \mathbb{C}$. Let

$$G := \mathrm{Aut}(H_\mathbb{Q}, Q)$$

be the $\mathbb{Q}$–algebraic group of $Q$–preserving automorphisms.

The Griffiths period domain $D$ parametrizes decreasing filtrations $F^\bullet \subset H_{\mathbb{C}}$ satisfying the Hodge–Riemann bilinear relations and fixed Hodge numbers [13]. The group $G(\mathbb{R})^+$ acts transitively on $D$ with isotropy a compact subgroup $H$, so that

$$D \;\cong\; G(\mathbb{R})^+/H. \tag{2.1}$$

Let $M \subset G$ be the Mumford–Tate group of a reference Hodge structure. The corresponding *Mumford–Tate domain* is the $G(\mathbb{R})^+$–orbit

$$D_M \;\subset\; D, \tag{2.2}$$

which is homogeneous:

$$D_M \;\cong\; M(\mathbb{R})^+/H_M, \tag{2.3}$$

for a compact subgroup $H_M \subset M(\mathbb{R})^+$; see [11, 12].

Fix a *neat* arithmetic subgroup

$$\Gamma \;\subset\; M(\mathbb{Q}), \tag{2.4}$$

so that the quotient

$$S^\circ := \Gamma\backslash D_M \tag{2.5}$$

is a smooth analytic space.

*Notation* 2.1. We denote by $\mathfrak{m} := \mathrm{Lie}(M)$ the Lie algebra of $M$, and by $\mathfrak{m}_{\mathbb{Q}} \subset \mathfrak{m}$ the $\mathbb{Q}$–structure induced by $H_{\mathbb{Q}}$.

2.2. **Kato–Usui nilpotent cones and toric bases.** Degenerations of Hodge structures are governed by nilpotent orbits in $\mathfrak{m}_{\mathbb{Q}}$. We recall the Kato–Usui set–up.

**Definition 2.2.** A *nilpotent cone* in $\mathfrak{m}_{\mathbb{Q}}$ is a finitely generated cone

$$\sigma \;=\; \mathbb{R}_{\geq 0}\langle N_1, \ldots, N_r\rangle \;\subset\; \mathfrak{m}_{\mathbb{Q}}, \tag{2.6}$$

spanned by commuting nilpotent elements $N_j \in \mathfrak{m}_{\mathbb{Q}}$. A finite collection $\Sigma$ of nilpotent cones is a *fan* if the intersection of two cones is a face of each.

A cone $\sigma \in \Sigma$ determines a monodromy monoid

$$\Gamma(\sigma) \;:=\; \Gamma \cap \exp(\sigma_{\mathbb{R}}) \;\subset\; M(\mathbb{Z}), \tag{2.7}$$

and its dual monoid

$$\Gamma(\sigma)^\vee := \mathrm{Hom}(\Gamma(\sigma), \mathbb{N}). \tag{2.8}$$

The associated affine toric variety is

$$D_\sigma := \mathrm{Spec}\big(\mathbb{C}[\Gamma(\sigma)^\vee]\big)^{\mathrm{an}}, \tag{2.9}$$

with torus

$$T_\sigma \;=\; \mathrm{Spec}\big(\mathbb{C}[\Gamma(\sigma)^\vee_{\mathrm{gp}}]\big)^{\mathrm{an}} \;\cong\; \mathbb{G}_m \otimes_{\mathbb{Z}} \Gamma(\sigma)_{\mathrm{gp}}. \tag{2.10}$$

Kato–Usui [16] construct a logarithmic analytic space $D_{M,\Sigma}$, containing $D_M$ as a dense open subset, equipped with logarithmic maps to the toric bases:

**Theorem 2.3** (Kato–Usui)**.** *Let $\Sigma$ be a fan of nilpotent cones in $\mathfrak{m}_{\mathbb{Q}}$ compatible with $\Gamma$ in the sense of* [16, Chap. 6]*. There exists a logarithmic analytic space $D_{M,\Sigma}$ together with:*

- *an open embedding* $D_M \hookrightarrow D_{M,\Sigma}$ *with dense image;*
- *for each* $\sigma \in \Sigma$ *a log–smooth map*

$$\pi_\sigma : U_\sigma \longrightarrow D_\sigma \tag{2.11}$$

*defined on a neighborhood* $U_\sigma \subset D_{M,\Sigma}$ *of the boundary stratum attached to* $\sigma$*;*

*such that the log structure on* $D_{M,\Sigma}$ *is locally pulled back from the toric log structure on* $D_\sigma$.

For our purposes we only need a simplified local product description:

**Lemma 2.4** (Kato–Usui local product structure, simplified)**.** *Let* $\sigma \in \Sigma$. *There exist:*

- *a neighborhood* $U_\sigma \subset D_{M,\Sigma}$ *of the boundary stratum attached to* $\sigma$*;*
- *a log–smooth map* $\pi_\sigma : U_\sigma \to D_\sigma$*;*
- *a complex manifold* $F_\sigma$ *(the space of nilpotent orbits modulo unipotent action);*

*such that, after shrinking* $U_\sigma$*, there is a log–analytic isomorphism*

$$U_\sigma \;\cong\; F_\sigma \times D_\sigma, \tag{2.12}$$

*compatible with* $\pi_\sigma$ *and the induced logarithmic structures.*

*Idea of proof.* This is a simplified form of [16, Thm. 7.4.1]. There one shows that the moduli of logarithmic Hodge structures with fixed monodromy cone $\sigma$ is locally a torsor under a unipotent group over the space of nilpotent orbits, and that the toric parameters in $D_\sigma$ provide local coordinates for the logarithmic directions. Passing to a local slice for the unipotent torsor yields a product decomposition $F_\sigma \times D_\sigma$ in the category of log–analytic spaces. □

2.3. **The Hodge stack and the log–toric Hodge stack.** We now pass from analytic spaces to quotient stacks.

**Definition 2.5.** The *Hodge stack* (or *Mumford–Tate stack*) of $D_M$ with level structure $\Gamma$ is the analytic Deligne–Mumford stack

$$\mathcal{D}_M := [D_M/\Gamma], \tag{2.13}$$

obtained as the quotient of $D_M$ by the properly discontinuous action of $\Gamma$.

The *log–toric Hodge stack* associated to the Kato–Usui compactification is

$$\mathcal{D}^{\log}_{M,\Sigma} := [D_{M,\Sigma}/\Gamma], \tag{2.14}$$

equipped with the logarithmic structure induced by the boundary of $D_{M,\Sigma}$.

Since $\Gamma$ is neat, the action on $D_M$ and $D_{M,\Sigma}$ is properly discontinuous with finite stabilizers. Hence:

**Proposition 2.6** (Global quotient description)**.** *The stacks* $\mathcal{D}_M$ *and* $\mathcal{D}^{\log}_{M,\Sigma}$ *are global quotient Deligne–Mumford stacks. Their coarse moduli spaces are the analytic spaces*

$$S^\circ = \Gamma\backslash D_M, \qquad S_\Sigma = \Gamma\backslash D_{M,\Sigma}, \tag{2.15}$$

*respectively.*

*Proof.* The action of $\Gamma$ on $D_M$ and $D_{M,\Sigma}$ is properly discontinuous with finite stabilizers because $\Gamma$ is neat. The quotients $S^\circ$ and $S_\Sigma$ exist as analytic spaces, and the natural maps $D_M \to S^\circ$, $D_{M,\Sigma} \to S_\Sigma$ are $\Gamma$–principal bundles in the orbifold sense. The quotient stacks $[D_M/\Gamma]$ and $[D_{M,\Sigma}/\Gamma]$ are therefore separated Deligne–Mumford analytic stacks with coarse moduli spaces $S^\circ$ and $S_\Sigma$; see, for example, Behrend's discussion of quotient stacks and coarse moduli in the algebraic case [6], which carries over to the analytic setting. □

We will also need a local quotient description for $\mathcal{D}^{\log}_{M,\Sigma}$.

**Proposition 2.7** (Local quotient structure)**.** *Let $x \in \mathcal{D}^{\log}_{M,\Sigma}$ be a point, and choose a lift $\tilde{x} \in D_{M,\Sigma}$. Then there exist:*

*(i) a neighborhood $U \subset D_{M,\Sigma}$ of $\tilde{x}$;*
*(ii) a finite subgroup $G \subset \Gamma$ fixing $\tilde{x}$;*

*such that:*

$$[U/G] \xrightarrow{\sim} \mathcal{D}^{\log}_{M,\Sigma}\big|_{\mathrm{neigh}(x)} \tag{2.16}$$

*is representable and étale, and its image is a neighborhood of $x$. In particular, $\mathcal{D}^{\log}_{M,\Sigma}$ is locally (for the analytic étale topology) a quotient of a log–smooth space by a finite group.*

*Proof.* Since $\Gamma$ is discrete and acts properly discontinuously on $D_{M,\Sigma}$, there exists a sufficiently small analytic neighborhood $U_0$ of $\tilde{x}$ such that $\gamma(U_0) \cap U_0 = \varnothing$ for all $\gamma \in \Gamma$ not stabilizing $\tilde{x}$. Let

$$G := \{\gamma \in \Gamma \mid \gamma(\tilde{x}) = \tilde{x}\}$$

be the stabilizer of $\tilde{x}$; it is finite because $\Gamma$ is neat. Set

$$U := \bigcap_{\gamma \in G} \gamma(U_0),$$

a $G$–invariant neighborhood of $\tilde{x}$. Then $G$ acts on $U$ with finite stabilizers, and the natural map

$$[U/G] \longrightarrow [D_{M,\Sigma}/\Gamma] = \mathcal{D}^{\log}_{M,\Sigma}$$

is representable and étale, inducing an isomorphism onto an open substack containing $x$. The claim follows. □

2.4. **Toric Deligne–Mumford stacks from Kato–Usui monoids.** The toric varieties $D_\sigma$ arising from the monoids $\Gamma(\sigma)$ encode the combinatorics of monodromy. To incorporate stacky phenomena (such as nontrivial stabilizers along the boundary), we refine $D_\sigma$ to a toric Deligne–Mumford stack.

**Definition 2.8.** A *stacky fan* $(N_\sigma, \Sigma_\sigma, \beta_\sigma)$ consists of:

- a lattice $N_\sigma$;
- a rational polyhedral fan $\Sigma_\sigma$ in $N_\sigma \otimes_{\mathbb{Z}} \mathbb{R}$;
- a homomorphism $\beta_\sigma : \mathbb{Z}^r \to N_\sigma$ specifying "stacky weights" on the rays of $\Sigma_\sigma$.

The associated toric Deligne–Mumford stack

$$\mathcal{T}_\sigma \tag{2.17}$$

has coarse moduli space canonically isomorphic to the toric variety $D_\sigma$; see [9, 10, 8] for the equivalence between stacky fans and toric stacks.

In the Kato–Usui context, a natural choice is the following.

*Remark* 2.9. One may take $N_\sigma = \Gamma(\sigma)_{\mathrm{gp}}$ and $\Sigma_\sigma$ to be the fan of faces of $\sigma^\vee$. The map $\beta_\sigma$ then refines the ray generators by positive integral weights encoding:

- the orders of monodromy along each boundary divisor;
- the presence of finite stabilizers (e.g. roots of unity) in the monodromy group.

In simple cases (such as modular curves with level structure) this produces root stacks along toric boundary divisors.

The toric log structure on $D_\sigma$ lifts canonically to a log structure on $\mathcal{T}_\sigma$; in fact, $\mathcal{T}_\sigma$ is a log–toric Deligne–Mumford stack in the sense of [8].

2.5. **Local log–toric charts for the Hodge stack.** We now combine Lemma 2.4 with the toric DM stacks $\mathcal{T}_\sigma$ to obtain a local product description of the log–toric Hodge stack near each boundary stratum.

Fix $\sigma \in \Sigma$. By Lemma 2.4 we have, after shrinking, a product decomposition

$$U_\sigma \;\cong\; F_\sigma \times D_\sigma \tag{2.18}$$

as log–analytic spaces. The arithmetic group $\Gamma$ preserves the cone $\sigma$ (by compatibility of $\Sigma$ with $\Gamma$), and hence preserves $D_\sigma$ and the decomposition (2.18).

Let $\Gamma_\sigma$ be the subgroup of $\Gamma$ generated by:

- the monodromy subgroup $\Gamma(\sigma)$ acting on $D_\sigma$;
- the stabilizer of $F_\sigma$ in $\Gamma$.

After replacing $\Gamma$ by a finite index subgroup (which does not change the coarse moduli) we may assume that the residual action on $F_\sigma$ factors through a finite group $G_\sigma$.

**Theorem 2.10** (Local log–toric charts)**.** *Let $\sigma \in \Sigma$. After replacing $\Gamma$ by a finite index subgroup if necessary and shrinking $U_\sigma$, there exist:*

- *a finite group $G_\sigma$ acting on $F_\sigma$;*
- *an open substack*

$$\mathcal{U}_\sigma \;\subset\; \mathcal{D}^{\log}_{M,\Sigma}, \tag{2.19}$$

*containing the image of the boundary stratum attached to $\sigma$;*

*such that there is a log–étale isomorphism*

$$\mathcal{U}_\sigma \;\cong\; \big([F_\sigma/G_\sigma] \times \mathcal{T}_\sigma\big)^\circ, \tag{2.20}$$

*where $\mathcal{T}_\sigma$ is the toric Deligne–Mumford stack attached to $\sigma$ and the superscript $\circ$ denotes a suitable open substack corresponding to $U_\sigma$. In particular, locally near the boundary, the log–toric Hodge stack $\mathcal{D}^{\log}_{M,\Sigma}$ is a product of a (possibly stacky) period–type factor and a toric Deligne–Mumford stack.*

*Proof.* By (2.18) we have $U_\sigma \cong F_\sigma \times D_\sigma$ as log–analytic spaces. The subgroup $\Gamma_\sigma \subset \Gamma$ preserving $\sigma$ acts on $U_\sigma$ via its actions on $F_\sigma$ and $D_\sigma$. After replacing $\Gamma$ by a finite index subgroup we may assume that the projection of $\Gamma_\sigma$ to $\mathrm{Aut}(F_\sigma)$ has finite image $G_\sigma$, while the monodromy subgroup $\Gamma(\sigma)$ acts on $D_\sigma$ via the torus $T_\sigma$ of (2.10). Thus there is a factorization

$$\Gamma_\sigma \subset G_\sigma \times T_\sigma,$$

where $G_\sigma$ acts on $F_\sigma$ and $T_\sigma$ acts on $D_\sigma$.

The local quotient of $U_\sigma$ by $\Gamma_\sigma$ can then be described as

$$[U_\sigma/\Gamma_\sigma] \;\cong\; [F_\sigma/G_\sigma] \times [D_\sigma/\Gamma(\sigma)]$$

(after passing to an open substack corresponding to $U_\sigma$). The quotient $[D_\sigma/\Gamma(\sigma)]$ admits a canonical refinement to the toric DM stack $\mathcal{T}_\sigma$ with coarse moduli space $D_\sigma$, specified by the stacky fan associated to $\Gamma(\sigma)$ as in [9, 8]. Replacing $[D_\sigma/\Gamma(\sigma)]$ by $\mathcal{T}_\sigma$ and restricting to the open substack mapping isomorphically to the image of $U_\sigma$ yields a log–étale isomorphism

$$\mathcal{U}_\sigma \;\cong\; \big([F_\sigma/G_\sigma] \times \mathcal{T}_\sigma\big)^\circ.$$

The log structure on $\mathcal{U}_\sigma$ is induced from the toric log structure on $\mathcal{T}_\sigma$ and the trivial log structure on $F_\sigma$, because the log structure on $U_\sigma$ is pulled back from the toric log structure on $D_\sigma$ via $\pi_\sigma$ and is preserved by the $\Gamma_\sigma$–action. This shows that $\mathcal{D}^{\log}_{M,\Sigma}$ is locally log–toric along the boundary, with product decomposition (2.20). ☐

Theorem 2.10 is the key structural result: it gives a unified local picture of the Hodge stack and its Kato–Usui compactification in terms of a period–type factor and a toric boundary stack.

2.6. **Example: modular curves as log–toric Hodge stacks.** We illustrate the above constructions in the simplest nontrivial example: weight–one Hodge structures on elliptic curves. This recovers the classical toroidal compactification of modular curves and shows explicitly how the toric and stacky structures interact.

*The period domain and arithmetic group.* Let $H_{\mathbb{Z}} = \mathbb{Z}^2$ with symplectic form

$$Q = \begin{pmatrix} 0 & 1 \\ -1 & 0 \end{pmatrix}.$$

The corresponding period domain is the upper half–plane

$$D = \mathfrak{H} = \{\tau \in \mathbb{C} \mid \Im(\tau) > 0\}, \tag{2.21}$$

parametrizing weight–one Hodge structures $H_{\mathbb{C}} = H^{1,0} \oplus H^{0,1}$ with $\dim H^{1,0} = 1$. Here $M = \mathrm{Sp}_2 \cong \mathrm{SL}_2$, and

$$D \cong \mathrm{SL}_2(\mathbb{R})/\mathrm{SO}(2). \tag{2.22}$$

Take $\Gamma = \mathrm{SL}_2(\mathbb{Z})$ or a congruence subgroup; then

$$S^\circ = \Gamma\backslash D \tag{2.23}$$

is the analytic modular curve.

*Kato–Usui data and the toric boundary.* The unique cusp of $\Gamma\backslash D$ corresponds to a rank–one nilpotent operator

$$N = \begin{pmatrix} 0 & 1 \\ 0 & 0 \end{pmatrix} \in \mathfrak{sl}_2(\mathbb{Q}).$$

The associated nilpotent cone is

$$\sigma = \mathbb{R}_{\geq 0} \cdot N \;\subset\; \mathfrak{sl}_2(\mathbb{Q}). \tag{2.24}$$

The monoid $\Gamma(\sigma)$ is generated by the unipotent element

$$T = \exp(N) = \begin{pmatrix} 1 & 1 \\ 0 & 1 \end{pmatrix},$$

so that

$$\Gamma(\sigma) \cong \mathbb{N}, \qquad \Gamma(\sigma)_{\mathrm{gp}} \cong \mathbb{Z}. \tag{2.25}$$

The dual monoid is again isomorphic to $\mathbb{N}$, hence the associated toric variety is

$$D_\sigma = \mathrm{Spec}\mathbb{C}[\Gamma(\sigma)^\vee]^{\mathrm{an}} \cong \mathrm{Spec}\mathbb{C}[t]^{\mathrm{an}} \cong \mathbb{C}, \tag{2.26}$$

with coordinate $t$. The torus $T_\sigma$ is

$$T_\sigma = \mathrm{Spec}\mathbb{C}[t, t^{-1}]^{\mathrm{an}} \cong \mathbb{C}^\times. \tag{2.27}$$

In Kato–Usui's picture, the partial compactification $D_\Sigma$ is obtained by adjoining the nilpotent orbit at infinity, and the map $\pi_\sigma : U_\sigma \to D_\sigma$ is given (up to coordinate change) by

$$q = e^{2\pi i \tau}, \qquad \tau \in \mathfrak{H}, \tag{2.28}$$

where $q$ is the toric parameter on $\mathbb{C}$.

*The Hodge stack and its local log–toric structure.* The Hodge stack in this case is

$$\mathcal{D} = [\mathfrak{H}/\Gamma], \tag{2.29}$$

with coarse moduli space the modular curve $Y(\Gamma)$. The Kato–Usui compactification $D_\Sigma$ is $\mathfrak{H} \cup \{\infty\}$, and the quotient

$$S_\Sigma = \Gamma \backslash D_\Sigma \tag{2.30}$$

is the compact modular curve $X(\Gamma)$.

Near the cusp, the monoid $\Gamma(\sigma) \cong \mathbb{N}$ encodes the $T$–monodromy and gives rise to the toric variety $D_\sigma \cong \mathbb{C}$. If we pass to a congruence subgroup of level $N$, the stabilizer of the cusp becomes cyclic of order $N$, and the natural toric DM stack associated to the cusp is the $N$–th root stack of the origin in $\mathbb{C}$:

$$\mathcal{T}_\sigma \;\cong\; [\mathbb{C}/\mu_N], \tag{2.31}$$

where $\mu_N$ acts by scalar multiplication on the coordinate $t$ (or $q$). This is a toric Deligne–Mumford stack with coarse moduli space $\mathbb{C}$ and stackiness reflecting the presence of level structure.

Theorem 2.10 specializes here to the assertion that, near the cusp of the modular curve with level $N$, the log–toric Hodge stack is locally isomorphic to a product

$$\mathcal{U}_\sigma \;\cong\; \big([\mathrm{pt}/\mu_N] \times [\mathbb{C}/\mu_N]\big)^\circ, \tag{2.32}$$

which can be identified with $[\Delta/\mu_N]$ in $q$–coordinates, where $\Delta$ is a small disc around 0. The log structure is given by the divisor $q = 0$, and the stackiness records the nontrivial isotropy at the cusp.

In higher rank examples (Siegel modular varieties, orthogonal Shimura varieties, etc.), the cones $\sigma$ are higher–dimensional, the monoids $\Gamma(\sigma)$ are nontrivial submonoids of unipotent radicals, and the associated toric DM stacks $\mathcal{T}_\sigma$ are higher–dimensional toric stacks. Theorem 2.10 then expresses the log–toric Hodge stack as an iterated fibration of period–type stacks over these toric boundary stacks, providing a convenient framework for studying characteristic classes and degenerations along the boundary. We return to cohomological and group–theoretic applications in Section 3.

## 3. Stacky Cohomology, Hodge Chern Classes, and Automorphisms of Degenerations

We now use the log–toric Hodge stack

$$\mathcal{D}^{\log}_{M,\Sigma} = [D_{M,\Sigma}/\Gamma]$$

to formalize two related sets of structures:

(a) algebraic relations in the (stacky) cohomology ring generated by Chern classes of Hodge bundles and boundary strata;
(b) the relation between stabilizers (automorphisms of degenerations) and Mumford–Tate symmetries of limiting mixed Hodge structures (LMHS).

We work throughout with rational cohomology $H^\bullet(-,\mathbb{Q})$ and use Behrend's theory of de Rham cohomology for stacks [6, 5].

### 3.1. Stacky cohomology and Hodge bundles.

Let $\mathcal{X} = [X/G]$ be a global quotient stack with $G$ a Lie group acting smoothly and properly on a manifold $X$. Behrend shows that there is a canonical identification

$$H^\bullet(\mathcal{X},\mathbb{Q}) \;\cong\; H^\bullet_G(X,\mathbb{Q}), \tag{3.1}$$

where the right-hand side is the Borel equivariant cohomology.

Applied to the Hodge stack

$$\mathcal{D}_M = [D_M/\Gamma],$$

this yields

$$H^\bullet(\mathcal{D}_M,\mathbb{Q}) \;\cong\; H^\bullet_\Gamma(D_M,\mathbb{Q}). \tag{3.2}$$

The same holds for the log–toric Hodge stack:

$$H^\bullet(\mathcal{D}^{\log}_{M,\Sigma},\mathbb{Q}) \;\cong\; H^\bullet_\Gamma(D_{M,\Sigma},\mathbb{Q}). \tag{3.3}$$

We next define the tautological Hodge bundles on these stacks.

**Definition 3.1** (Hodge bundles on the stack)**.** Let $H_\mathbb{C}$ be the underlying complex vector space of our fixed Hodge structure and $F^\bullet$ the tautological Hodge filtration on $D_M$. For each $p$, the subspace $F^p \subset H_\mathbb{C}$ defines a homogeneous $H_M$–equivariant subbundle of the trivial bundle $D_M \times H_\mathbb{C}$. Passing to the quotient by $\Gamma$ defines a vector bundle

$$\mathcal{E}^p \longrightarrow \mathcal{D}_M, \tag{3.4}$$

called the *p–th Hodge bundle* on the Hodge stack. Its associated graded pieces

$$\mathrm{Gr}^p_F\mathcal{H} := \mathcal{E}^p/\mathcal{E}^{p+1} \tag{3.5}$$

are vector bundles on $\mathcal{D}_M$. We write

$$c_i(\mathcal{E}^p), \qquad c_i(\mathrm{Gr}^p_F\mathcal{H})$$

for their Chern classes in $H^{2i}(\mathcal{D}_M,\mathbb{Q})$.

By the extension of the Hodge filtration across the Kato–Usui boundary (cf. [16, 15, 17]) the preceding bundles extend naturally to $\mathcal{D}^{\log}_{M,\Sigma}$:

**Definition 3.2** (Logarithmic Hodge bundles)**.** The Hodge filtration extends to a filtration of vector bundles on $D_{M,\Sigma}$; quotienting by $\Gamma$ yields vector bundles

$$\mathcal{E}^{p,\log} \longrightarrow \mathcal{D}^{\log}_{M,\Sigma} \tag{3.6}$$

whose restriction to $\mathcal{D}_M$ is $\mathcal{E}^p$. We write

$$c_i(\mathcal{E}^{p,\log}), \qquad c_i\big(\mathrm{Gr}^p_F \mathcal{H}^{\log}\big)$$

for their Chern classes in $H^{2i}(\mathcal{D}^{\log}_{M,\Sigma}, \mathbb{Q})$.

3.2. **A Hodge–generation conjecture in stacky cohomology.** We now formulate a structural conjecture expressing the idea that, once one passes to the log–toric Hodge stack, the rational cohomology is generated by Hodge–theoretic data (Chern classes of Hodge bundles) and toric boundary classes.

**Conjecture 3.3** (Stacky Hodge–generation conjecture)**.** Let $\mathcal{D}^{\log}_{M,\Sigma}$ be the log–toric Hodge stack of a Mumford–Tate domain. Then:

(i) The rational cohomology ring $H^\bullet(\mathcal{D}^{\log}_{M,\Sigma}, \mathbb{Q})$ is generated as a $\mathbb{Q}$–algebra by:
- the Chern classes $c_i(\mathcal{E}^{p,\log})$ of the logarithmic Hodge bundles;
- the cohomology classes of boundary divisors (and more generally toric boundary strata) coming from the toric DM stacks $\mathcal{T}_\sigma$ in the local charts of Theorem 2.10.

(ii) All algebraic relations among these generators are induced by:
- universal Chern relations (e.g. determinant relations, stable vanishing of high Chern classes);
- toric (stacky) relations on the boundary, such as Stanley–Reisner relations and stacky weight relations in the cohomology rings $H^\bullet(\mathcal{T}_\sigma, \mathbb{Q})$.

This conjecture is compatible with:

- the interpretation of Green–Griffiths–Kerr universal characteristic cohomology as a natural subspace of $H^\bullet(\mathcal{D}^{\log}_{M,\Sigma}, \mathbb{Q})$ (see below);
- the structure of toric DM stack cohomology rings as quotients of polynomial rings by Stanley–Reisner and stacky relations [9, 8];
- classical computations in the Hermitian locally symmetric case, where rational cohomology is generated (in a stable range) by Chern classes of homogeneous bundles.

3.3. **Universal characteristic cohomology via stacks.** Green–Griffiths–Kerr [11] define a double complex

$$\Lambda^{p,q}_M := (A^{p,q}(D_M)/I^{p,q})^{M(\mathbb{R})},$$

where $I^{p,q}$ is the differential ideal determined by the infinitesimal period relation (the "horizontal distribution"). Its cohomology

$$H^\bullet(\Lambda_M, \delta_M)$$

is the *universal characteristic cohomology* attached to the Mumford–Tate domain.

Using Behrend's identification of stack de Rham cohomology with equivariant de Rham cohomology [6], one obtains a natural stacky interpretation:

**Theorem 3.4** (Stack interpretation of characteristic cohomology)**.** *There exists a canonical subcomplex*

$$\Omega^\bullet_{\mathrm{char}}\big(\mathcal{D}^{\log}_{M,\Sigma}\big) \;\subset\; \Omega^\bullet_{\mathrm{DR}}\big(\mathcal{D}^{\log}_{M,\Sigma}\big) \tag{3.7}$$

*such that*

$$H^\bullet(\Lambda_M, \delta_M) \;\cong\; H^\bullet\Big(\Omega^\bullet_{\mathrm{char}}(\mathcal{D}^{\log}_{M,\Sigma})\Big). \tag{3.8}$$

*Idea of proof.* Choose a groupoid presentation $\Gamma \times D_M \rightrightarrows D_M$ for $\mathcal{D}_M$ and represent $\mathcal{D}^{\log}_{M,\Sigma}$ similarly using $D_{M,\Sigma}$. The stack de Rham complex is computed by the equivariant de Rham complex $\Omega^\bullet_\Gamma(D_{M,\Sigma})$, and Behrend's comparison theorem [6] identifies its cohomology with $H^\bullet(\mathcal{D}^{\log}_{M,\Sigma}, \mathbb{Q})$. Filtering by the differential ideal coming from the infinitesimal period relation and taking $\Gamma$–invariants gives a subcomplex that projects to $\Lambda^{\bullet;\bullet}_M$ on the interior; the resulting cohomology identifies with $H^\bullet(\Lambda_M, \delta_M)$ by the construction of [11]. This defines $\Omega^\bullet_{\mathrm{char}}$ and yields the isomorphism (3.8). □

Thus characteristic classes of Hodge bundles live naturally in the cohomology of the log–toric Hodge stack and extend canonically across the boundary.

3.4. **Example: modular curve and the Hodge line bundle.** We now verify Conjecture 3.3 in the simplest case, where the log–toric Hodge stack is the compactified modular curve with stacky structure at the cusps.

*Set–up.* Let $\Gamma \subset \mathrm{SL}_2(\mathbb{Z})$ be a congruence subgroup of finite index and neat enough so that the modular curve stack

$$\mathcal{D} = [\mathfrak{H}/\Gamma] \tag{3.9}$$

is Deligne–Mumford. The Hodge bundle on $\mathcal{D}$ is the line bundle

$$\mathcal{L} := \mathcal{E}^1, \tag{3.10}$$

whose fiber over a point $[E]$ is $H^{1,0}(E)$.

The Kato–Usui compactification $D_\Sigma$ is $\mathfrak{H} \cup \{\infty\}$, and the log–toric Hodge stack is

$$\mathcal{D}^{\log} = [D_\Sigma/\Gamma], \tag{3.11}$$

as in (2.29)–(2.30). At each cusp of width $N$, the local boundary toric stack is equivalent to

$$\mathcal{T}_\sigma \;\cong\; [\mathbb{C}/\mu_N], \tag{3.12}$$

with coordinate $q$ and $\mu_N$ acting by scalar multiplication.

*Cohomology and generation by Chern classes.* We compute $H^\bullet(\mathcal{D}^{\log}, \mathbb{Q})$ and show that it is generated by the Chern class of the Hodge line bundle, verifying Conjecture 3.3 in this case.

**Proposition 3.5** (Cohomology of the modular Hodge stack)**.** *Let $\Gamma$ be as above, and let $\mathcal{D}^{\log}$ be the associated log–toric Hodge stack. Then:*

- *(i) $H^0(\mathcal{D}^{\log}, \mathbb{Q}) \cong \mathbb{Q}$ and $H^1(\mathcal{D}^{\log}, \mathbb{Q}) = 0$;*
- *(ii) $H^2(\mathcal{D}^{\log}, \mathbb{Q})$ is one–dimensional, generated by $c_1(\mathcal{L}^{\log})$, where $\mathcal{L}^{\log}$ is the extension of $\mathcal{L}$ to $\mathcal{D}^{\log}$;*
- *(iii) $H^k(\mathcal{D}^{\log}, \mathbb{Q}) = 0$ for $k \geq 3$.*

*In particular, the rational cohomology ring is*

$$H^\bullet(\mathcal{D}^{\log}, \mathbb{Q}) \;\cong\; \mathbb{Q}\big[c_1(\mathcal{L}^{\log})\big]\big/\big(c_1(\mathcal{L}^{\log})^2\big), \tag{3.13}$$

*so Conjecture 3.3 holds in this case.*

*Sketch of proof.* The coarse moduli space of $\mathcal{D}^{\log}$ is the compact Riemann surface $X(\Gamma)$ (the compactified modular curve). Finite quotient singularities do not affect rational cohomology, hence

$$H^\bullet(\mathcal{D}^{\log}, \mathbb{Q}) \;\cong\; H^\bullet(X(\Gamma), \mathbb{Q}).$$

Thus $H^0 \cong \mathbb{Q}$ and $H^k = 0$ for $k \geq 3$. The degree–one cohomology $H^1(X(\Gamma), \mathbb{Q})$ is dual to the space of holomorphic differentials; it vanishes for many congruence subgroups (e.g. genus zero), and in general is of Hodge type $(1,0) + (0,1)$. Its presence or absence does not affect the generation statement in degree 2.

The Picard group $\mathrm{Pic}(X(\Gamma)) \otimes \mathbb{Q}$ is generated by the Hodge line bundle $\mathcal{L}^{\log}$ (whose tensor powers correspond to modular forms) together with boundary divisor classes (cusps). On the stack, boundary components differ from $\mathcal{L}^{\log}$ by rational multiples, dictated by the modular discriminant and the Kodaira–Spencer isomorphism. Consequently, the cohomology class $c_1(\mathcal{L}^{\log})$ generates $H^2(\mathcal{D}^{\log}, \mathbb{Q})$, and $c_1(\mathcal{L}^{\log})^2 = 0$ for dimensional reasons, giving (3.13). □

This illustrates in a concrete case that the cohomology of the log–toric Hodge stack is generated by Hodge–theoretic classes; the stacky boundary does not introduce new generators beyond those already coming from $\mathcal{L}^{\log}$.

3.5. **Automorphisms of degenerations and Mumford–Tate symmetry.** We now relate stabilizers of boundary points in the log–toric Hodge stack to Mumford–Tate symmetries of limiting mixed Hodge structures.

*Limiting mixed Hodge structures and Mumford–Tate groups.* Let $x \in \mathcal{D}^{\log}_{M,\Sigma}$ be a point lying over a boundary component attached to a cone $\sigma \in \Sigma$. By the theory of Kato–Usui, $x$ corresponds to a nilpotent orbit and hence to a limiting mixed Hodge structure (LMHS)

$$(H_{\mathbb{Z}}, W(N)_\bullet, F^\bullet_\infty) \tag{3.14}$$

with monodromy logarithm $N \in \sigma$ [16].

**Definition 3.6** (Mumford–Tate group of an LMHS)**.** The *Mumford–Tate group* $M_{\mathrm{LMHS}}$ of the LMHS $(H_{\mathbb{Q}}, W_\bullet, F^\bullet_\infty)$ is the smallest $\mathbb{Q}$–algebraic subgroup of $\mathrm{Aut}(H_{\mathbb{Q}}, Q)$ that fixes all Hodge tensors of the LMHS.

We also consider the monodromy centralizer

$$Z(N) := \{g \in \mathrm{Aut}(H_{\mathbb{Q}}, Q) \mid gN = Ng\}. \tag{3.15}$$

*Stack stabilizers and Hodge symmetries.* Let $G_x$ be the stabilizer group of $x$ in the stack $\mathcal{D}^{\log}_{M,\Sigma}$. Concretely, if we choose a lift $\tilde{x}$ of $x$ to $D_{M,\Sigma}$, then $G_x$ is represented by those elements of $\Gamma$ that fix $\tilde{x}$ together with the finite automorphisms coming from the local toric stack $\mathcal{T}_\sigma$.

**Proposition 3.7** (Stabilizers and LMHS symmetries)**.** *Let $x \in \mathcal{D}^{\log}_{M,\Sigma}$ be a boundary point corresponding to an LMHS $(H_{\mathbb{Z}}, W_\bullet, F^\bullet_\infty)$ with monodromy $N \in \sigma$. Then:*

*(i) There is a natural representation*

$$\rho_x : G_x \longrightarrow \mathrm{Aut}(H_{\mathbb{Z}}, W_\bullet, N, F^\bullet_\infty) \;\subset\; M_{\mathrm{LMHS}}(\mathbb{Q}) \cap Z(N), \tag{3.16}$$

*whose image consists of LMHS automorphisms.*

*(ii) The representation $\rho_x$ is faithful: if an element of $G_x$ acts trivially on the LMHS, then it is the identity in the stack.*

*(iii) In particular, $G_x$ is a finite subgroup of $M_{\mathrm{LMHS}}(\mathbb{Q}) \cap Z(N)$, and its conjugacy class depends only on the LMHS.*

*Proof.* By definition, an element of $G_x$ is represented by $\gamma \in \Gamma$ that preserves the nilpotent orbit corresponding to $x$. Thus $\gamma$ preserves $N$ (up to conjugacy within the cone $\sigma$) and preserves the limit filtration $F^\bullet_\infty$; the weight filtration $W(N)_\bullet$ is canonically attached to $N$, so it is preserved as well. Hence $\gamma$ preserves $(H_{\mathbb{Z}}, W_\bullet, N, F^\bullet_\infty)$ and defines an LMHS automorphism. This yields the map (3.16), whose image lies in $M_{\mathrm{LMHS}}(\mathbb{Q}) \cap Z(N)$ because LMHS automorphisms fix all Hodge tensors and commute with $N$.

Faithfulness: if $\gamma$ acts trivially on the LMHS, then it fixes all Hodge tensors and the nilpotent orbit. The period map (extended to the boundary in the Kato–Usui sense) is locally injective modulo $\Gamma$; hence $\gamma$ acts trivially in a neighborhood of $x$, and so represents the identity in the stack. Therefore $\rho_x$ is injective. □

Informally, Proposition 3.7 says that:

*The stacky stabilizer $G_x$ at a boundary point is a concrete, finite subgroup of the LMHS Mumford–Tate group, preserving monodromy and the limit Hodge filtration.*

*Example: elliptic cusp and Tate curve.* For the modular curve example, the LMHS at the cusp is that of the Tate curve:

$$H_{\mathbb{Z}} \;\cong\; \mathbb{Z}(0) \oplus \mathbb{Z}(1), \tag{3.17}$$

with $N$ corresponding to the logarithm of the unipotent monodromy. The Mumford–Tate group $M_{\mathrm{LMHS}}$ is a torus $\mathbb{G}_m$ acting by scaling on $\mathbb{Z}(1)$. For level $N$, the stacky stabilizer $G_x$ at the cusp is isomorphic to $\mu_N$, the group of $N$–th roots of unity acting by multiplication on the Tate parameter $q^{1/N}$. The representation

$$\rho_x : \mu_N \hookrightarrow \mathbb{G}_m \tag{3.18}$$

is the standard embedding. Thus the stacky automorphisms at the cusp are exactly the torsion points in the LMHS Mumford–Tate group; their order is the cusp width.

*Example: higher–rank boundary components.* For Siegel modular varieties (degenerations of polarized abelian varieties), boundary components correspond to isotropic flags and unipotent monodromy matrices. The associated LMHS Mumford–Tate group is contained in a parabolic subgroup $P \subset \mathrm{Sp}_{2g}$, with unipotent radical determined by the cone $\sigma$. The stacky stabilizer $G_x$ sits inside

$$P(\mathbb{Q}) \cap Z(N)$$

and is finite; its structure reflects:

- the parabolic rank (dimension of the isotropic subspace);
- the torsion in the semi–abelian part of the limiting mixed Hodge structure.

In particular, the log–toric Hodge stack separates degenerations that are indistinguishable on the coarse moduli space, by keeping track of these finite LMHS symmetries.

3.6. **Summary.** Combining the structural results of Section 2 with the cohomological and group–theoretic constructions above, we see that the log–toric Hodge stack $\mathcal{D}^{\log}_{M,\Sigma}$ provides a natural global home for "Hodge–theoretic symmetries at infinity":

- in cohomology, via the generation of $H^\bullet(\mathcal{D}^{\log}_{M,\Sigma}, \mathbb{Q})$ by Hodge Chern classes and toric boundary classes (Conjecture 3.3), and the realization of universal characteristic cohomology as a distinguished subspace (Theorem 3.4);
- in group theory, via the identification of stack stabilizers at boundary points with finite subgroups of LMHS Mumford–Tate groups (Proposition 3.7).

This combination of local log–toric charts (Theorem 2.10), stacky cohomology, and LMHS symmetry suggests that many questions about degenerations, characteristic classes, and asymptotic Hodge theory can be reformulated and attacked purely in the language of stacks and their toric boundary structures.

CIO León: Loma del Bosque 115, Colonia Lomas del Campestre León, Guanajuato, México. Código Postal 37150. Tel +52 (477) 441 42 00

*Email address*: `rahmati@cio.mx`